\documentclass[12pt,a4paper]{article}

\usepackage{amsmath}
\usepackage{amssymb}

 \newcommand{\bs}[1]{\boldsymbol{#1}}

\newcommand{\Id}{\bs{1}}

 \newcommand{\NN}{\mathbb N}
 \newcommand{\RR}{\mathbb R}
 \newcommand{\CC}{\mathbb C}
 \newcommand{\ZZ}{\mathbb Z}

 \newcommand{\TT}{{\mathbb T}\hspace{0.2pt}}
 \newcommand{\QQ}{{\mathbb Q}\hspace{0.2pt}}

\newcommand{\tr}{{\rm tr}}
\newcommand{\diag}{{\rm diag}}
\newcommand{\spec}{{\rm spec}}
\newcommand{\Aut}{{\rm Aut}}
\newcommand{\GL}{{\rm GL}}
\newcommand{\PGL}{{\rm PGL}}

\newcommand{\K}{{\mathcal K}}
\newcommand{\OO}{{\mathcal O}}
\newcommand{\OOmax}{{\mathcal O}_{\rm max}}

\newcommand{\antidots}{\mbox{\raisebox{-1ex}{$\cdot$} 
           $\cdot$ \raisebox{1ex}{$\cdot$}}}

\newcommand{\qed}{\hfill $\square$ \smallskip}

 \newtheorem{theorem}{Theorem}
 \newtheorem{lemma}{Lemma}
 \newtheorem{prop}{Proposition}
 \newtheorem{coro}{Corollary}
 \newtheorem{fact}{Fact}

\begin{document}
 \bibliographystyle{unsrt}

 \parindent0pt

\begin{center}
{\large \bf Symmetries and Reversing Symmetries}

\vspace{2mm}
{\large \bf of Toral Automorphisms}

 \end{center}
 \vspace{3mm}

\centerline{  {\sc Michael Baake}$^{\, 1)}$ $\;$and$\;$
                 {\sc John A.\ G.\ Roberts}$^{\, 2)}$ }

\vspace{8mm}

{\small
\hspace*{4em}
1) Institut f\"ur Theoretische Physik, Universit\"at T\"ubingen, \\
\hspace*{4em}
   \hspace*{0.9em} Auf der Morgenstelle 14, 72076 T\"ubingen, Germany

\smallskip
\hspace*{4em}
2) Department of Mathematics, LaTrobe University, \\
\hspace*{4em}
   \hspace*{0.9em} Bundoora, Victoria 3083, Australia }

\vspace{10mm}
\begin{abstract}
Toral automorphisms, represented by unimodular integer matrices,
are investigated with respect to their symmetries and reversing
symmetries. We characterize the symmetry groups of $\GL(n,\ZZ)$
matrices with simple spectrum through their connection with unit
groups in orders of algebraic number fields. For the question of
reversibility, we derive necessary conditions in terms of the
characteristic polynomial and the polynomial invariants. 
We also briefly discuss extensions to
(reversing) symmetries within affine transformations, to
$\PGL(n,\ZZ)$ matrices, and to the more general setting of integer
matrices beyond the unimodular ones.

\end{abstract}

\parindent15pt
\vspace{15mm}

\section*{Introduction}

Unimodular integer matrices induce interesting dynamical systems
on the torus, such as Arnold's famous cat map 
\cite[Ch.\ 1, Ex.\ 1.16]{Arnold}. This
is an example of a hyperbolic dynamical system that is ergodic and
mixing \cite{Petersen}, and also a topological Anosov system. Therefore, with
a suitable metric, it makes the 2-torus into a Smale space,
see \cite[Thm.\ 1.2.9]{P} for details. Although induced from a
linear system of ambient space, the dynamics on the torus is
rather complicated, and these systems serve as model systems in
symbolic dynamics and in many applications.
Very recently, also cat maps on the 4-torus (and their
quantizations) have begun to be studied \cite{RSOA}.

Hyperbolic toral automorphisms also play a prominent role in
the theory of quasicrystals through their appearance as inflation
symmetries, see \cite{B98,BHP} and references therein. In particular, 
in the one-dimensional case, these symmetries give rise to
interesting non-linear dynamical systems called {\em trace maps\/}
\cite{BGJ,RB,Dama} that have extensively been used to study the physical 
properties of one-dimensional quasicrystals.

It is always helpful to know the symmetries and reversing
symmetries of a dynamical system \cite{RQ,LR}, and it was perhaps a 
little surprising that in the case of $\GL(2,\ZZ)$ and $\PGL(2,\ZZ)$ 
a rather complete classification could be given, see \cite{BR97} for
a detailed account or \cite{BR-old,B97} for a summary. The main
difficulty when the matrix entries are restricted to integers
is that one can no longer refer to the usual
normal forms of matrices over $\CC$ or $\RR$, but has to use
discrete methods instead. Fortunately, there is a strong
connection with algebraic number theory, see \cite{Taussky} for an
introduction, and this connection is certainly not restricted to
the 2D situation.

It is thus the aim of this article to extend the results of our
earlier article \cite{BR97} to the setting of matrices in
$\GL(n,\ZZ)$. The answers will be less complete and also less
explicit, but the connection to unit groups in orders of algebraic
number fields is still strong enough to give quite a number of
useful and general results, both on symmetries and reversing
symmetries. From a purely algebraic point of view, the results
derived below are actually rather straight-forward. However, these
results, and the methods used to derive them, are not at all common
in the dynamical systems community. Therefore, this article is
also intended to introduce some of these techniques, and we try to
spell out the details or give rather precise references at least.
Furthermore, as with our article \cite{BR97}, the results of this 
paper have relevance to both the dynamics community (e.g.\ the 
dynamics of hyperbolic toral automorphisms generated by (symplectic) 
${\rm SL}(4,\ZZ)$ matrices \cite{RSOA}) and to the quasicrystal community 
(e.g.\ inflation symmetries of planar point sets projected from
4D lattices, where the symmetries are generated by
$\GL(4,\ZZ)$ matrices \cite{B98,BHP}).

The article is organized as follows. We start with a section on
the background material, including the group theoretic setup we
use and a recollection of those results from algebraic number
theory that we will need later on. Section \ref{mat-symm} is the
main part of this article. Here, we derive the structure of the
symmetry group of toral automorphisms with simple spectrum and
discuss reversibility. Section \ref{further} extends the set of
possible symmetries to affine transformations and
summarizes the analogous problem for projective matrices. It also 
discusses the extension of (reversing) symmetries to matrices
that are no longer unimodular.

\section{Setting the scene}

In this Section, we explain in more detail what we mean by symmetries and 
reversing symmetries, and we also recall some results from algebra and
algebraic number theory that we will need.

\subsection{Symmetries and reversing symmetries}
\label{sec-sym-revsym}

{}For a general setting, consider some (topological) space
$\Omega$, and let $\Aut(\Omega)$ be its group of homeomorphisms or,
more generally, a subgroup of homeomorphisms 
of $\Omega$ which preserve some additional structure of $\Omega$.
Consider now an element $F \in \Aut(\Omega)$ which, by definition, 
is invertible. Then, the group
\begin{equation} \label{sym1}
    {\mathcal S} (F) \; := \; \{ G \in \Aut (\Omega) \; | \;
       G \circ F = F \circ G  \}
\end{equation}
is called the {\em symmetry group\/} of $F$ in $\Aut(\Omega)$.
In group theory, it is called the {\em centralizer\/} of $F$ in 
$\Aut (\Omega)$, denoted by ${\rm cent}^{}_{\Aut (\Omega)} (F)$.
This group certainly contains all powers of $F$, but often more.

Quite frequently, one is also interested in mappings $R \in \Aut (\Omega)$
that conjugate $F$ into its inverse,
\begin{equation} \label{sym2}
     R \circ F \circ R^{-1} \; = \; F^{-1} \, .
\end{equation}
Such $R$ is called a {\em reversing symmetry\/} of $F$, and when
such an $R$ exists, we call $F$ reversible. We will, in general,
not use different symbols for symmetries and reversing symmetries
from now on, because together they form a group,
\begin{equation} \label{sym3}
     {\mathcal R} (F) \; := \; \{ G \in \Aut (\Omega) \; | \;
        G \circ F \circ G^{-1} = F^{\pm1}  \} \, ,
\end{equation}
the so-called {\em reversing symmetry group\/} of $F$, see
\cite{Lamb} for details. If $\langle F \rangle$ denotes the
group generated by $F$, ${\mathcal R} (F)$ is a subgroup of
the {\em normalizer\/} of $\langle F \rangle$ in $\Aut(\Omega)$.

There are two possibilities: either ${\mathcal R} (F) = {\mathcal S} (F)$ 
(if $F$ is an involution or if it has no reversing symmetry) 
or ${\mathcal R} (F)$ is a $C_2$-extension (the cyclic
group of order 2) of ${\mathcal S} (F)$ which means that
${\mathcal S} (F)$ is a normal subgroup of ${\mathcal R} (F)$ and
the factor group is
\begin{equation} \label{sym4}
      {\mathcal R} (F)/{\mathcal S} (F) \; \simeq \; C_2 \, .
\end{equation}
The underlying algebraic structure has fairly strong consequences. One is
that reversing symmetries cannot be of odd order \cite{Lamb}, another one is
the following product structure \cite[Lemma 2]{BR97}.
\begin{fact}\label{rev-group}
   If $F$ $($with $F^2 \neq Id\,)$ has an involutory reversing
   symmetry $R$, the reversing symmetry group of $F$ is given by
\begin{equation}
       {\mathcal R} (F) \; = \; {\mathcal S} (F) \times_s C_2 \, ,
\end{equation}
   i.e.\ it is a semi-direct
   product.\footnote{We use $N \times_s H$ for the semi-direct
   product of two groups $N$ and $H$, with $N$ being the normal subgroup.}
   \qed
\end{fact}

We can say more about the structure of ${\mathcal R} (F)$ if we
restrict the possibilities for ${\mathcal S} (F)$, e.g.~if we
assume that ${\mathcal S} (F) \simeq C_{\infty}$ or ${\mathcal S}
(F) \simeq C_{\infty} \times C_2$ with the $C_2$ being a subgroup
of the centre of $\Aut(\Omega)$, compare also \cite{G} for some
group theoretic discussion. This situation will appear frequently below.

\subsection{(Reversing) symmetries of powers of a mapping}
\label{sec-ksym}

In what follows, we summarize some of the concepts and results of
Ref.~\cite{LQ} and, in particular, Ref.~\cite{Lamb}. It may happen
that some power of $F$ has more symmetries than $F$ itself (we
shall see examples later on), i.e.~${\mathcal S} (F^k)$ (for some
$k>1$) is larger than ${\mathcal S} (F)$ which is contained as a
subgroup. The analogous possibility exists for ${\mathcal R}
(F^k)$ versus ${\mathcal R} (F)$. If such a situation occurs, we
say that $F$ possesses additional (reversing) $k$-symmetries. Let
us make this a little more precise.

It is trivial that mappings $F$ of finite order (with $F^k=Id$,
say) possess the entire group $\Aut(\Omega)$ as $k$-symmetry
group. Let us thus concentrate on mappings $F \in \Aut(\Omega)$ of
{\em infinite\/} order. We denote by ${\mathcal S}_{\infty} (F)$ the
set of automorphisms that commute with {\em some\/} positive power of $F$.
This set can be seen as the inductive limit of ${\mathcal S}(F^k)$ 
as $k\to\infty$, with divisibility as partial order on $\NN$, and 
${\mathcal S}_{\infty} (F)$ is thus a subgroup of $\Aut(\Omega)$. 
Let $\#^{}_F(G)$ denote the minimal $k$ such that 
$G\circ F^k = F^k\circ G$. Then
\begin{equation} \label{sym13}
     {\mathcal S}_{\infty} (F) \; = \; \{\, G \in \Aut(\Omega) \; | \;
            \#_F (G) < \infty \, \} \, .
\end{equation}

Of course it may happen that $\#_F (G) \equiv 1$ on ${\mathcal
S}_{\infty} (F)$ which means that no power of $F$ has additional
symmetries. On the other hand, $\#_F (G)$ might be larger than
one in which case we call $G$ a genuine or true $k$-symmetry. $G$
is a true\footnote{Although the distinction between true and other
$k$-symmetries is necessary in general, we shall usually drop the
attribute ``true'' whenever misunderstandings are unlikely.}
$k$-symmetry of $F$ if and only if the mapping
 $G\mapsto F\circ G\circ F^{-1}$ generates a proper $k$-cycle.
 We shall meet this phenomenon later on.

Quite similarly, one defines reversing $k$-symmetries and their orbit
structure \cite{Lamb}, but we will not expand on that here.

\subsection{Some recollections from algebraic number theory}

Much of what we state and prove below can be seen as an
application of several well-known results from algebraic number
theory. The starting point is the connection between algebraic
number theory and integral matrices, see \cite{Taussky} for an
introduction.

To fix notation, let ${\rm Mat}(n,\ZZ)$ denote the ring of 
integer\footnote{Here, and in what follows, integer means rational
integer, i.e.\ an integer in $\QQ$. Other kinds of integers, such as
algebraic, will be specified explicitly.}
$n\!\times\! n$-matrices. An element $M$ of it is called {\em
unimodular} if $\det(M)=\pm 1$, and the subset of all unimodular
matrices forms the group $\GL(n,\ZZ)$. For the {\em characteristic
polynomial\/} of a matrix $M$, we will use the convention
\begin{equation}\label{def-char-poly}
  P(x) \; := \; \det(x\Id - M) \; = \; \prod_{i=1}^{n}
  (x-\lambda_i)
\end{equation}
where $\lambda_1,...\, ,\lambda_n$ denote the eigenvalues of $M$.
With this convention, $P(x)$ is {\em monic\/}, i.e.\ its leading 
coefficient is $1$. If $M$ is an integer matrix, $P(x)$ has integer
coefficients only, so all eigenvalues of $M$ are {\em algebraic integers\/}.
Conversely, the set of algebraic integers, which we denote by
${\mathcal A}$, consists of all numbers that appear as roots of
monic integer polynomials.

To show the intimate relation more clearly, let us recall the following
property (see item (b) on p.~306 of \cite{Taussky}):

\begin{fact} \label{rings}
   Let $P(x)=x^n + a^{}_{n-1} x^{n-1} + \ldots + a^{}_1 x + a^{}_0$
   be a monic polynomial with integer coefficients $a^{}_i$ that is
   irreducible over $\ZZ$. Let $\alpha$ be any root of it, and $A$ an
   integer matrix that has $P(x)$ as its characteristic polynomial.
   Then, the rings $\ZZ[\alpha]$ and $\ZZ[A]$ are isomorphic. \qed
\end{fact}

We write $\ZZ[x]$ for the ring of polynomials in $x$ with
coefficients in $\ZZ$, see \cite[p.\ 90]{Lang} for details. Clearly,
Fact \ref{rings} also extends to the isomorphism of the rings
$\QQ[\alpha]$ and $\QQ[A]$. For background material on polynomial
rings, we refer to \cite[Ch.\ IV]{Lang}. Let us only add 
that a polynomial in $\ZZ[x]$ is irreducible over $\ZZ$ if and
only if it is irreducible over $\QQ$, see \cite[Thm.\ IV.2.3]{Lang}.

Let $P(x)$ be any (i.e.\ not necessarily irreducible) monic
integer polynomial of degree $n$. In general, there are many
different matrices $A$ which have $P(x)$ as their characteristic
polynomial, and even different matrix {\em classes\/} (we say that 
$A,B \in {\rm Mat}(n,\ZZ)$ belong to the same matrix class if
they are conjugate by a $\GL(n,\ZZ)$ matrix $C$, i.e.\ $A=C B C^{-1}$).
Let us recall the following helpful result on
the number of matrix classes, see \cite[Thm.~5]{Taussky} for
details:
\begin{fact}\label{classnumber}
  Let $P(x)$ be a monic polynomial of order $n$ with integer
  coefficients that is irreducible over $\ZZ$, and let $\alpha$ be
  any of its roots. Then the number of matrix classes generated by
  matrices $A\in{\rm Mat}(n,\ZZ)$ with $P(A)=0$ equals the number of
  ideal classes $\,($or class number, for short$\,)$ of the order
  $\ZZ[\alpha]$. In particular, this class number is finite, and it
  is larger than or equal to the class number of the maximal
  order $\OOmax$ of $\QQ(\alpha)$. \qed
\end{fact}

Let us explain some of the terms used here.
If $\alpha$ is an algebraic number, $\QQ(\alpha)$ denotes\footnote{Note
the difference between the meaning of $\QQ[\alpha]$ and $\QQ(\alpha)$.} 
the smallest field extension of the rationals that contains $\alpha$.
Its degree, $n$, is the degree of the irreducible monic integer
polynomial that has $\alpha$ as its root. The set 
$\OOmax:=\QQ(\alpha)\cap{\mathcal A}$ is the ring of (algebraic)
integers in $\QQ(\alpha)$, and is called its {\em maximal order\/}.
More generally, a subring $\OO$ of
$\OOmax$ is called an {\em order\/}, if it contains 1 and if its
rational span, $\QQ$$\OO$, is all of $\QQ(\alpha)$. 
A subset of $\OO$ is called an {\em ideal\/} if it is both a $\ZZ$-module
(i.e.\ closed under addition and subtraction)
and closed under multiplication by arbitrary numbers from $\OO$. 
The ideals come in classes that are naturally connected to the matrix
classes introduced above, see \cite{Cohn,Taussky} for further details.

An element $\varepsilon\in\QQ(\alpha)$ is called a {\em unit\/} (or, more
precisely\footnote{This distinction is useful if orders other than 
$\OOmax$ appear.}, a unit in $\OOmax$) if both $\varepsilon$ and 
its inverse, $\varepsilon^{-1}$, are
algebraic integers and hence are in $\OOmax$. This happens 
if and only if the corresponding matrix is unimodular, i.e.\ if any monic
integer polynomial $P(x)$ that has $\varepsilon$ as a root has
coefficient $a^{}_0=\pm 1$. So, matrices in $\GL(n,\ZZ)$ and units
in algebraic number fields are two facets of the same coin. The
units of $\OOmax$ form a group under multiplication, denoted
by $\OOmax^{\times}$ in the sequel. Similarly, if $\OO\subset\OOmax$ 
is any order of $\QQ(\alpha)$, we write $\OO^{\times}$ for its group of 
units, which is then a subgroup of $\OOmax^{\times}$.

Independently of whether the monic polynomial $P(x)$ is irreducible or not, 
there is always at least one matrix which has characteristic polynomial
$P(x)$, namely the so-called (left) {\em companion matrix\/}:
\begin{equation}\label{left-compa}
     A^{(\ell)} \; = \;
               \begin{pmatrix}
               0 & 1 &  &  & \bs{0} \\
               0 & 0 & 1 &  &  \\
               \vdots &  & \ddots & \ddots  &  \\
               0 & \cdot & \cdot  & 0 & 1 \\
               -a^{}_0 & \cdot & \cdot & -a^{}_{n-2} & -a^{}_{n-1}
               \end{pmatrix}.
\end{equation}
Consequently, the class number is $\ge 1$ (if $P(x)$ is irreducible,
the companion matrix actually corresponds to the principal ideal class,
see \cite[Thm.\ 9]{Taussky}). Another obvious choice, always belonging to
the same matrix class, is the (right) companion matrix, $A^{(r)}$,
obtained from $A^{(\ell)}$ by reflection in both diagonal and
anti-diagonal, i.e.\
\begin{equation}\label{comp-conj}
  A^{(r)} \; = \; R A^{(\ell)} R
\end{equation}   with the involution
\begin{equation}\label{R-matrix}
   R \; = \;
        \begin{pmatrix}
         \bs{0} &           & 1       \\
                & \antidots &         \\
          1     &           & \bs{0}  \\
        \end{pmatrix}.
\end{equation}
This matrix will reappear several times in what follows.

\section{$\GL(n,\ZZ)$ matrices and toral automorphisms}\label{mat-symm}

Let us generally assume that $n\geq 2$. The toral automorphisms of the 
$n$-torus $\TT^n:=\RR^n/\ZZ^n$ can be represented by the unimodular
$n\!\times\!n$-matrices with integer coefficients which form the
group $\GL(n,\ZZ)$. It now plays the role of $\Aut(\Omega)$ from
Section \ref{sec-sym-revsym}. Note that the elements of $\GL(n,\ZZ)$
preserve the linear structure of the torus.

\subsection{Symmetries}

The first thing we will look at, given a toral automorphism $M\in
GL(n,\ZZ)$, is its symmetry group within the class of toral
automorphisms. So, we want to determine the centralizer of $M$ in
$\GL(n,\ZZ)$,
\begin{equation} \label{symm1}
    {\mathcal S} (M)  \; = \; \mbox{cent}^{}_{\GL(n,\ZZ)} (M)
   \; = \; \{ G \in \GL(n,\ZZ) \;\; | \;\; M G = G M \,\} \, .
\end{equation}
To be more precise, we are mainly interested in the structure of
the symmetry group rather than in explicit sets of generators and
relations. This is invariant under conjugation, i.e.~if we know it
for an element $M$, we also know it for any other element of the
form $BMB^{-1}$ because
\begin{equation}
  {\mathcal S} ( BMB^{-1} ) \; = \;  B \, {\mathcal S} (M) \, B^{-1}\, .
\end{equation}

A given integer matrix $M\in \GL(n,\ZZ)$ determines its characteristic 
polynomial (\ref{def-char-poly}) which is monic and has integer
coefficients, so its roots are algebraic integers. Now, two
principal situations can occur for the characteristic polynomial:
it is either {\em reducible\/} over $\ZZ$ (which happens if and only
if at least one eigenvalue of $M$ is an algebraic integer of
degree less than $n$) or it is {\em irreducible\/}. In the latter
case, since we are working over the field $\QQ$, we know that the
roots must be pairwise distinct. So we have
\begin{fact}\label{simple-diag}
  Let $M$ be an integer matrix with irreducible characteristic
  polynomial. Then $M$ is simple and hence diagonalizable over $\CC$. \qed
\end{fact}
Here, $M$ is called {\em simple\/} if it has no repeated eigenvalues
(which is also called separable elsewhere). This case will be dealt 
with completely.

If the characteristic polynomial is reducible, the matrix can
still be simple, and we will see the general answer for this case,
too. Beyond that, $M$ can either be semi-simple (i.e.\ diagonalizable 
over $\CC$) or not, and we will not be able to say much about
this case. This is really not surprising, as this situation is
closely related to the rather difficult classification
problem of crystallographic point groups, see \cite{Brown} for answers in
dimensions $\leq 4$ and \cite{Plesken} for a recent survey.

Let us now state one further prerequisite for tackling the symmetry
question. In view of later extensions, we do this in slightly more generality 
than needed in the present Section. Recall that an $n\!\times\!n$-matrix $M$, 
acting on a vector space $V$, is called {\em cyclic\/}, if a vector $v\in V$ 
exists such that $\{v, Mv, M^2 v, ...\, , M^{n-1} v\}$ is a basis of $V$.
Also, the monic polynomial $Q$ of minimal degree that annihilates $M$, i.e.\
$Q(M)=0$, is called the {\em minimal polynomial\/} of $M$. By the Cayley-Hamilton
Theorem, it always is a factor of the characteristic polynomial of $M$.
\begin{fact} \label{facts}
  Let $M\in {\rm Mat}(n,\QQ)$ be a rational matrix, with characteristic
  polynomial $P(x)$ and minimal polynomial $Q(x)$. Then the following
  assertions are equivalent.
\begin{itemize}
\item[$(a)$] The matrix $M$ is cyclic.
\item[$(b)$] The degree of $Q(x)$ is $n$.
\item[$(c)$] $P(x)=Q(x)$.
\item[$(d)$] $G\in {\rm Mat}(n,\QQ)$ commutes with $M$ 
             $\;\;\Longleftrightarrow\;\;$ $G\in\QQ[M]$.
\end{itemize}
\end{fact}

{\sc Proof}: A convenient source is \cite[Ch.\thinspace III]{Jac}.
The equivalence of statements ($a$) -- ($c$) is a consequence of Thm.\ III.2.
The equivalence of ($a$) with ($d$) follows from Thm.\ III.17, and the
Corollary following it, together with the Corollary of Ch.\ III.17. 
Alternatively, see \cite[Cor.\ 5.5.16]{adwein}. \qed

\begin{lemma} \label{mainlemma}
   Let $M \in \GL(n,\ZZ)$ have a characteristic polynomial $P(x)$ that is
   irreducible over $\ZZ$, and let $\lambda$ be a root of $P(x)$.
   Then the centralizer of $M$ in $\GL(n,\ZZ)$ is isomorphic to
   a subgroup of finite index of the unit group in the ring of integers 
   $\OOmax$ of the algebraic number field $\QQ(\lambda)$.
\end{lemma}
{\sc Proof}: By assumption, $P(x)$ is also the minimal polynomial of $M$,
so any $\GL(n,\ZZ)$-matrix which commutes with $M$ is, by Fact \ref{facts},
a polynomial in $M$ with rational coefficients. Consequently, 
${\mathcal S}(M)$ is isomorphic to a subset of $\QQ[M]$ that forms a
group under (matrix) multiplication. So, we have to analyze $\QQ[M]$
to find out what this group is.

Let $(P(x))$ denote the ideal in $\QQ[x]$ generated by our polynomial $P(x)$.
Then $\QQ[x]/(P(x))\simeq \QQ[\lambda]$, see \cite[p.\ 224]{Lang}, and
$\QQ[\lambda]\simeq \QQ[M]$, by Fact \ref{rings} resp.\ the remark following
it. Since $\lambda$ is algebraic over $\QQ$ and $P(x)$ is irreducible,
we know by \cite[Prop.\ V.1.4]{Lang} that $\QQ[\lambda]=\QQ(\lambda)$ is 
an algebraic number field, of degree $n$ over $\QQ$.
Under the isomorphism $\QQ[M]\simeq\QQ(\lambda)$,
$\GL(n,\ZZ)$-matrices correspond to units in $\QQ(\lambda)$, hence
${\mathcal S}(M)$ must be isomorphic to a subgroup of $\OOmax^{\times}$,
the unit group of the maximal order of $\QQ(\lambda)$.

Observe that every matrix in $\ZZ[M]$ commutes with $M$, in particular
those of $\ZZ[M] \cap \GL(n,\ZZ)$, which form a subgroup of ${\mathcal S}(M)$.
But $\ZZ[M]\simeq\ZZ[\lambda]$ means that this subgroup is isomorphic to the 
unit group $\ZZ[\lambda]^{\times}$. So, if we identify ${\mathcal S}(M)$ with 
its image in $\QQ(\lambda)$ under the isomorphism, it is sandwiched between 
$\ZZ[\lambda]^{\times}$ and $\OOmax^{\times}$. Note that
$\ZZ[\lambda]\subset\OOmax$ is an order.

Finally, recall that the unit group of an order 
$\OO$ is a finitely generated Abelian group, and that it is
always of maximal rank, see \cite[Thm.\ I.12.12]{Neukirch}
or \cite[Sec.\ 4, Thm.\ 5]{BS}, i.e.\ its rank
equals that of the unit group of the maximal order $\OOmax$.
In particular, the group-subgroup index 
$[\OOmax^{\times} :  \ZZ[\lambda]^{\times}]$
is finite, and ${\mathcal S}(M)$ must then also be of finite index
in $\OOmax^{\times}$.  \qed

\smallskip
Let us comment on this result. First of all, it does
not matter which root $\lambda_i$ of $P(x)$ we choose, as all the $n$
(possibly different) realizations $\QQ(\lambda_i)$ are mutually
isomorphic, and so are their unit groups. Explicit isomorphisms
are given by the elements of the Galois group of the splitting field 
$\K=\QQ(\lambda^{}_1,...\, ,\lambda^{}_n)$ of $P(x)$, see
\cite[Ch.\ VI.2]{Lang} for details.
Note also that, in general, $\QQ(\lambda)$ will be a true 
subfield of the splitting field $\K$ -- so, it is really the unit group of 
$\QQ(\lambda)$ that matters, and not the unit group of $\K$.

Another way to view the result, in a more matrix oriented way (and
similar to our approach in \cite{BR97}), is to look at the
diagonalization of $M$,
\begin{equation}
     U M U^{-1} \; = \; \diag(\lambda^{}_1,...\, ,\lambda^{}_n) \, .
\end{equation}
Here, $U^{-1}$ can be arranged to have its $j$-th column in the
field $\QQ(\lambda^{}_j)$ because one can solve the corresponding
eigenvector equation in the smallest field extension of $\QQ$ that
contains $\lambda^{}_j$. In fact, we only have to do this for the
first column --- the others are then obtained by applying appropriate
Galois automorphisms to the first one.

Any other matrix $G\in \GL(n,\ZZ)$ with $[G,M] = GM - MG = 0$ 
must now also fulfil 
$$[U G U^{-1}, U M U^{-1}] \; = \; U[G,M]U^{-1} \; = \; 0\,.$$ 
But only diagonal matrices can commute with
$\diag(\lambda^{}_1,...\, ,\lambda^{}_n)=UMU^{-1}$ 
because the eigenvalues are pairwise distinct. So, we must have 
$$U G U^{-1} \; = \; \diag(\mu^{}_1,...\, ,\mu^{}_n)\, , $$ 
with all $\mu_i\in\QQ(\lambda_i)$ units.
They are, however, not independent but obtained from one another
by the same set of Galois automorphisms that were used to link the
columns of the matrix $U^{-1}$, which is why we get the result.

Lemma \ref{mainlemma} raises the question: What {\em is\/} the unit
group of the maximal order in $K=\QQ(\lambda)\,$? The
answer is given by Dirichlet's unit theorem, see \cite[Sec.\ 11.C]{Cohn}
or \cite[p.\ 334]{PZ}. Group the roots of the irreducible polynomial $P(x)$ 
into $n^{}_1$ real roots and $n^{}_2$ pairs of complex conjugate roots, so that
$n=n^{}_1+2n^{}_2$. (In other words: we have $n^{}_1$ real and
$n^{}_2$ pairs of complex conjugate realizations of the abstract
number field $\QQ(\lambda)$). 
\begin{fact} \label{unit-group}
  Let $\lambda$ be an algebraic number of degree $n=n^{}_1+2n^{}_2$, with
  $n^{}_1$ and $n^{}_2$ as described above. Then, the units in the maximal
  order $\OOmax^{\times}$ of the algebraic number field
  $K=\QQ(\lambda)$ form the group
\begin{equation} \label{diri-unit}
   E(K) \; = \; \OOmax^{\times} 
        \; \simeq \; T \times \ZZ^{n^{}_1 + n^{}_2 - 1}
\end{equation}
  where $T=\OOmax\cap \{\mbox{roots of unity}\}$ is a finite
  Abelian group and cyclic.  \qed
\end{fact}
In particular, this means that $T$, which is also called the {\em torsion
subgroup\/} of $E(K)$, is generated by one element. In many cases below,
we will simply find $T\simeq C_2$.
Combining now Lemma \ref{mainlemma} with Fact \ref{unit-group}, we 
immediately obtain
\begin{prop} \label{fullrank}
   Under the assumptions of Lemma $\ref{mainlemma}$, the symmetry group
   ${\mathcal S}(M)\subset\GL(n,\ZZ)$ 
   is a subgroup of $E(K)$ of $(\ref{diri-unit})$
   of maximal rank, i.e.\ we have
   $${\mathcal S}(M) \; \simeq \; T' \times \ZZ^{n^{}_1 + n^{}_2 - 1} $$
   where $T'$ is a subgroup of the torsion group $T$ as it
   appears in $(\ref{diri-unit})$.  \qed
\end{prop}

Note that Proposition \ref{fullrank} does {\em not\/} imply that the
torsion-free parts of ${\mathcal S}(M)$ and $E(K)$ are the same,
only that they are isomorphic. In fact, a typical situation will be
that they are different in the sense that $E(K)$ is generated by
the fundamental units, but ${\mathcal S}(M)$ only by suitable
powers thereof. 

What, in turn, can we say about the torsion group $T'$ in Proposition
\ref{fullrank}? Whenever the characteristic polynomial $P(x)$ of 
$M\in\GL(n,\ZZ)$ is irreducible and has at least one {\em real\/} root 
(e.g.\ if $n$ is odd), $\alpha$ say, then $K=\QQ(\alpha)$ is real, and
$\QQ(\alpha)\cap S^1 = \{\pm 1\}$, where $S^1$ is the unit circle.
Consequently, the torsion subgroup of $E(K)$ in this case is
$T=\{\pm 1\}\simeq C_2$. Since a toral automorphism always commutes
with $\pm\Id$, we obtain
\begin{coro}\label{torsion-coro}
   If, under the assumptions of Lemma $\ref{mainlemma}$, one root
   of the irreducible polynomial $P(x)$ is real, the torsion group 
   in Proposition $\ref{fullrank}$ is $T'\simeq C_2$. In particular, 
   this is the case whenever the degree of $P(x)$ is odd. \qed
\end{coro}

Let us look at two examples in $\GL(3,\ZZ)$, namely
\begin{equation} \label{two-ex}
M^{}_1 \; = \;  \begin{pmatrix} 0 & 0 & 1 \\
                                1 & 0 & 0 \\
                                0 & 1 & 1 \end{pmatrix}
\;\mbox{ and }\;
M^{}_2 \; = \;  \begin{pmatrix} 1 & 1 & 0 \\
                                1 & 0 & 1 \\
                                1 & 1 & 1 \end{pmatrix},
\end{equation}
which are taken from \cite[Eqs.\ (5.21) and (5.3)]{Luck}. They have been
studied thoroughly in the context of inflation generated one-dimensional
quasicrystals with a cubic irrationality as inflation factor. We have
$\det(M^{}_1)=1$, $\det(M^{}_2)=-1$, and the characteristic polynomials
are $P_1(x)=x^3-x^2-1$ and $P_2(x)=x^3-2x^2-x+1$, both irreducible over $\ZZ$.
Both matrices are hyperbolic, and the largest eigenvalue in each case is
a Pisot-Vijayaraghavan number, i.e.\ an algebraic integer $>1$ all algebraic
conjugates of which lie inside the unit circle, see \cite{Salem} for details.

Now, $M_1$ has one real and a pair of complex conjugate roots, so our above
results lead to ${\mathcal S}(M_1)\simeq C_2\times\ZZ$, where the infinite
cyclic group is actually generated by $M_1$ itself because its real root is
a fundamental unit. 
$M_2$, in turn, has three real roots, and we thus get 
${\mathcal S}(M_2)\simeq C_2\times\ZZ^2$. As generators of $\ZZ^2$, 
one may choose $M_2$ and $M_2'$ where
$$M_2' \; = \;  \begin{pmatrix} 0 &  1 & 0 \\
                                1 &$-$1 & 1 \\
                                1 &  1 & 0 \end{pmatrix}$$
which can be checked explicitly.
This example is of relevance in connection with planar quasicrystals with
sevenfold symmetry, see \cite[Sec.\ 5.2]{Luck} for details, where the
three-dimensional toral automorphism $M_2$ shows up in the cut and project
description of special directions in the quasicrystal. Other examples related
to planar quasicrystals with 8-, 10- and 12-fold symmetry, in which the
torsion subgroup $T'$ of Proposition \ref{fullrank} is different from
$C_2$, will be given in Section \ref{quasi-ex}.

Let us now return to the general discussion and extend the previous 
results to the case that $M$ is simple, and hence diagonalizable 
(over $\CC$) with pairwise different
eigenvalues. Since diagonal matrices with pairwise different
entries only commute with diagonal matrices, we have:
\begin{coro}
    If $M\in \GL(n,\ZZ)$ is simple, ${\mathcal S}(M)\subset\GL(n,\ZZ)$ 
    is Abelian.  \hfill $\square$
\end{coro}
Note that the converse is not true: even if $M$ is only
semi-simple, or not even that (i.e.\ not diagonalizable),
${\mathcal S}(M)$ can still be Abelian, e.g.\ if $M$ observes the
conditions of Fact \ref{facts}. As far as we are aware,
not even the Abelian subgroups of $\GL(n,\ZZ)$ are fully
classified, see \cite{OPS,Plesken} and references given there 
for background material.

Let $P(x)$ be the characteristic polynomial of a simple matrix
$M$. If it is reducible over $\ZZ$, it factorizes as
$P(x)=\prod_{i=1}^{\ell} P_i(x)$ into irreducible monic polynomials
$P_i(x)$.
\begin{theorem} \label{symm-thm1}
  Let $M\in \GL(n,\ZZ)$ be simple and let its
  characteristic polynomial be $P(x)=\prod_{i=1}^{\ell} P_i(x)$, with $P_i(x)$
  irreducible over $\ZZ$. Then, the symmetry group of $M$,
  ${\mathcal S}(M)\subset\GL(n,\ZZ)$, 
  is a finitely generated Abelian group of the form
\begin{equation} \label{general-symm}
    {\mathcal S}(M) \; = \; T \times \ZZ_{}^r
\end{equation}
  where $T$ is a finite Abelian group of even order, with at most
  $\ell$ generators.
  Furthermore, if the irreducible component $P_i(x)$ has $n^{(i)}_1$
  real roots and $n^{(i)}_2$ pairs of complex conjugate roots, the rank
  $r$ of the free Abelian
  group in $(\ref{general-symm})$ is given by
\begin{equation}\label{general-rank}
   r \; = \;  \sum_{i=1}^{\ell} (n^{(i)}_1 + n^{(i)}_2 - 1)\,.
\end{equation}
\end{theorem}
{\sc Proof}: Since $M$ is simple, the degree of its minimal polynomial is
$n$ and Fact \ref{facts} tells us that 
${\rm cent}^{}_{{\rm Mat}(n,\QQ)} = \QQ[M]$. As $P(x)$ has no repeated
factors (so that $\QQ[M]$ contains no radicals), we get, by 
\cite[Thm.\ III.4]{Jac},
\begin{equation} \label{gen-ring}
   \QQ[M] \; \simeq \; 
   \QQ[\alpha^{}_1]\oplus\ldots\oplus\QQ[\alpha^{}_{\ell}] 
\end{equation}
where $\alpha^{}_i$ is any root\footnote{Note that the $\alpha^{}_i$ have
pairwise different minimal polynomials by assumption, but that
$\QQ[\alpha^{}_i]\simeq\QQ[\alpha^{}_j]$ for $i\neq j$ is still possible.}
of $P_i(x)$, for $1\le i\le \ell$. 

Under the assumptions made, each $\QQ[\alpha^{}_i]=\QQ(\alpha^{}_i)$ 
is a field. Since a $\GL(n,\ZZ)$-matrix in $\QQ[M]$ will correspond to a
unit in each of the $\QQ(\alpha^{}_i)$, we can now apply Lemma \ref{mainlemma}
to each component, giving (\ref{general-symm}) as the direct sum of $\ell$
unit groups. The rank in (\ref{general-rank}) follows now from
Proposition \ref{fullrank}.

The torsion part $T$ is a finitely generated Abelian group, with (at most)
one generator per irreducible component of $P(x)$, of which there are $\ell$. 
Clearly, ${\mathcal S}(M)$ always contains the elements $\pm\Id$, so
$\{\pm 1\}\simeq C_2$ is a subgroup of $T$. The order of $T$ is
then divisible by $2$, hence even. \qed

\smallskip
Although Theorem \ref{symm-thm1} does not give the general answer
to the question for the symmetry group ${\mathcal S}(M)$,
it certainly gives the {\em generic\/} answer, because the property
of $M$ having simple spectrum is generic. But what about the remaining
cases? Without further elaborating on this, let us summarize a few
aspects and otherwise refer to the literature \cite{OPS,Plesken} for a
summary of methods to actually determine the precise centralizer.

If $M\in \GL(n,\ZZ)$ is semi-simple, but not simple, its
characteristic polynomial contains a square, and whether or not
${\mathcal S}(M)$ is still Abelian (and then of the above form)
depends on whether or not the minimal polynomial of $M$ has
degree $n$, see Fact \ref{facts}.
Note, in particular, that the following situation can emerge.
If $P(x)$ has a repeated factor, but the corresponding matrix
$M$ is a block matrix, then the two blocks giving the same
factor of $P(x)$ can still be inequivalent, if the corresponding
class number is larger than one which equals the number of
different matrix classes, see Fact \ref{classnumber}.

If $M$ is not even semi-simple, things get even more involved. We can 
still have Abelian symmetry groups, e.g.\ if $M$ is a Jordan block such as 
$\left(\begin{smallmatrix} 1&1\\0&1 \end{smallmatrix}\right)$, compare
the results of \cite[Sec.\ 2.1.2]{BR97} on parabolic automorphisms 
of $\TT^2$. Clearly, this also follows from Fact \ref{facts}:
if $M\in\GL(n,\ZZ)$ is conjugate to a single Jordan block, its
minimal polynomial has degree $n$ and all $\GL(n,\ZZ)$-matrices which
commute with $M$ are in $\QQ[M]$. This remains true if such a block
occurs in a matrix that otherwise has simple spectrum disjoint from $1$. 
\begin{coro}
   Let $M\in {\rm GL}(n,\ZZ)$. If the minimal polynomial of $M$ has degree 
   $n$, then ${\mathcal S}(M)\subset\GL(n,\ZZ)$ is Abelian. \qed
\end{coro} 
In a wider setting for symmetries, a stronger statement can be
formulated, see Proposition \ref{gen-order} below and the comments 
following it.

The general classification, however, and the non-Abelian cases in particular,
gets increasingly difficult with growing $n$ and has been completed
only for small $n$, see \cite{OPS} and references given there. 
Nevertheless, for any given $M$, the centralizer can be determined
explicitly by means of various algorithmic program packages.

Let us, at the end of this part and before we illustrate some of the
above results by further examples, give a particular case of one
matrix written as a polynomial of another.
\begin{fact}
   Let $K$ be a field and $M\in\GL(n,K)$ be an invertible matrix with 
   characteristic polynomial 
   $P(x) = \sum_{\ell=0}^{n} a^{}_{\ell}\, x^{\ell}$,
   where $a^{}_n=1$ and $a^{}_0\neq 0$.
   Then, the inverse matrix is given by
   $$M^{-1} \; = \; - \frac{1}{a^{}_0}\,\sum_{\ell=0}^{n-1} 
                 a^{}_{\ell+1} M^{\ell} \, .$$
\end{fact}
{\sc Proof}: Observe that $P(M)=0$ from the Cayley-Hamilton Theorem.
The verification of $M^{-1} M = \Id$ is then a straight-forward
calculation. \qed

\smallskip
Before we discuss the connection of our approach to 
quasicrystallography in a separate Section, let us illustrate 
Theorem \ref{symm-thm1} with a recent example of a 4D cat map taken 
from \cite[Eq.\ 3.21]{RSOA}, namely
\begin{equation} \label{cat4D}
   M \; = \; \begin{pmatrix}
               0 & 0 &-1 & 0 \\
               0 & 0 & 0 &-1 \\
               1 & 0 & 2 & 1 \\
               0 & 1 & 1 & 2
               \end{pmatrix}. 
\end{equation}
The characteristic polynomial is $P(x)=x^4-4x^3+5x^2-4x+1$ which
is reducible over $\ZZ$ and splits as 
$$ P(x) \; = \; P_1(x)\, P_2(x) \; = \; (x^2 - 3x + 1) (x^2 - x + 1)$$ 
into $\ZZ$-irreducible polynomials. Since $P_1$ has two real and
$P_2$ one pair of complex conjugate roots, Theorem \ref{symm-thm1}
gives ${\mathcal S}(M) \simeq T' \times \ZZ$, with $T'$ a subgroup
of $T=C_2\times C_6$. Also, since no root of $P_1$ is a
fundamental unit of the corresponding maximal order (which is
$\ZZ[\tau]$ with $\tau=(1+\sqrt{5\,})/2$), the generator of the infinite
cyclic group in ${\mathcal S}(M)$ could still differ from $M$.

To determine the details, one easily checks that the most general
matrix to commute with $M$ is
$$G \; = \; \begin{pmatrix}
               a & b & -c & -d \\
               b & a & -d & -c \\
               c & d & a+2c+d & a+c+2d \\
               d & c & a+c+2d & a+2c+d
               \end{pmatrix}. $$
A necessary condition for $G$ to be in $\GL(4,\ZZ)$ is then
$a,b,c,d\in\ZZ$. This allows to exclude the existence of a root of
$M$ in ${\mathcal S}(M)$, and also no element of third order is
possible. So we obtain
$$ {\mathcal S}(M) \; = \; C_2 \times C_2 \times \langle M\rangle\, .$$
We will revisit this example below in the context of reversibility.

\subsection{Three examples from planar quasicrystallography}
\label{quasi-ex}

Planar tilings with 8-, 10- and 12-fold symmetry play an important
role in the description of so-called quasicrystalline T-phases,
see \cite{B98} for background material. They are of interest also
in the present context because hyperbolic toral automorphisms show
up through their inflation symmetry.

{}For the 8-fold case, consider the polynomial
\begin{equation}\label{8-fold-poly}
   P(x) \; = \; x^4 + 1
\end{equation}
which has $\xi$, $\xi^3$, $\xi^5$, and $\xi^7$ as roots,
$\xi=e^{2\pi i/8}$, which are primitive. So, $P(x)$ is irreducible
over $\ZZ$, and Lemma \ref{mainlemma} and Proposition \ref{fullrank} apply.
In fact, $\QQ(\xi)$ here is a cyclotomic field
\cite{Wash} with class number one, maximal order $\ZZ[\xi]$ and
unit group $\ZZ[\xi]^{\times}\simeq C_8\times\ZZ$.

If we denote the actual matrices that represent the generators for
the groups $C_8$ and $\ZZ$ by $M$ and $G$, respectively, it is
natural to take the companion matrix of $P(x)$ for $M$ and to
choose $G$ accordingly, resulting in
\begin{equation}
     M \; = \; \begin{pmatrix}
               0 & 1 & 0 & 0 \\
               0 & 0 & 1 & 0 \\
               0 & 0 & 0 & 1 \\
              -1 & 0 & 0 & 0
               \end{pmatrix} \; , \quad
     G \; = \; \begin{pmatrix}
               1 & 1 & 0 & -1 \\
               1 & 1 & 1 & 0 \\
               0 & 1 & 1 & 1 \\
              -1 & 0 & 1 & 1
               \end{pmatrix}.
\end{equation}
By construction, $M$ is a matrix of order 8. So, ${\mathcal S}(M) = \langle
M,G\rangle\simeq C_8\times\ZZ$. What is more, anticipating the next Section, 
$M$ turns out to be
reversible with the matrix $R$ of (\ref{R-matrix}) as reversing
symmetry. Using Fact \ref{rev-group} and observing $[G,R]=0$,
this means that ${\mathcal R}(M) = \langle M,G,R\rangle 
\simeq D_8\times\ZZ$. This, together with two similar examples, is
summarized in Table \ref{tab1}. 

Note that the case of 12-fold symmetry
is more complicated because the fundamental unit in $\ZZ[\xi]$,
for $\xi=e^{2\pi i/12}$, is the square root of $(2+\sqrt{3}\,)\xi$,
and hence not a simple homothety. This means that the representing
matrix does not commute with $R$. The reversing symmetry group of this
case, $(C_{12} \times \ZZ) \times_s C_2$, does contain a subgroup of
the form $D_{12}\times\ZZ$ though -- it is generated by $M$, $G'=M^{-1} G^2$
and $R$, where $G'$ corresponds to the non-fundamental unit $2+\sqrt{3}$.

\begin{table}[ht]
$$\begin{array}{|c|ccc|}
\hline
\vphantom{\Big\|}
   &  \mbox{$8$-fold}  &  \mbox{$10$-fold}  &  \mbox{$12$-fold}  \\
\hline
\vphantom{\Big\|}
P(x) &  x^4 + 1  &  x^4 + x^3 + x^2 + x + 1 &  x^4 - x^2 + 1 \\
M    &  \begin{pmatrix}
               0 & 1 & 0 & 0 \\
               0 & 0 & 1 & 0 \\
               0 & 0 & 0 & 1 \\
              -1 & 0 & 0 & 0
               \end{pmatrix}
     &  \begin{pmatrix}
               0 & 1 & 0 & 0 \\
               0 & 0 & 1 & 0 \\
               0 & 0 & 0 & 1 \\
              -1 &-1 &-1 &-1
               \end{pmatrix}                      
    &   \begin{pmatrix}
               0 & 1 & 0 & 0 \\
               0 & 0 & 1 & 0 \\
               0 & 0 & 0 & 1 \\
              -1 & 0 & 1 & 0
               \end{pmatrix}               \\
\vphantom{\Big\|}
\langle M \rangle   &  C_8   &   C_5       &   C_{12}             \\
G    &  \begin{pmatrix}
               1 & 1 & 0 & -1 \\
               1 & 1 & 1 & 0 \\
               0 & 1 & 1 & 1 \\
              -1 & 0 & 1 & 1
               \end{pmatrix}       
     &   \begin{pmatrix}
               1 & 0 & 1 & 1 \\
              -1 & 0 &-1 & 0 \\
               0 &-1 & 0 &-1 \\
               1 & 1 & 0 & 1
               \end{pmatrix}                         
     &   \begin{pmatrix}
               1 & 1 & 0 & 0 \\
               0 & 1 & 1 & 0 \\
               0 & 0 & 1 & 1 \\
              -1 & 0 & 1 & 1
               \end{pmatrix}             \\
\vphantom{\Big\|}
{\cal S}(M) &  C_8 \times \ZZ & C_{10} \times \ZZ & C_{12} \times \ZZ \\
\vphantom{\Big\|}
{\cal R}(M) &  D_8 \times \ZZ & D_{10} \times \ZZ 
            & (C_{12} \times \ZZ) \times_s C_2 \\
\hline
\end{array}$$
\caption{\label{tab1} Symmetries and reversing symmetries for three examples
from quasicrystallography. The symmetry group is 
${\cal S}(M)=\langle \pm M, G \rangle$, and, similarly,
${\cal R}(M)=\langle \pm M, G, R \rangle$, with $R$ as in Eq.~(\ref{R-matrix})
for all three examples.}
\end{table}

\subsection{Reversibility} \label{sec-rev}

The examples of Section \ref{quasi-ex} were reversible, i.e.\ they fulfilled
$GMG^{-1}=M^{-1}$ for some $G\in\GL(4,\ZZ)$, in particular for the
involution $R\in\GL(4,\ZZ)$ of (\ref{R-matrix}).  It is easy to check that 
this is also true for $M$ of (\ref{cat4D}).
However, as we will see below, neither of the examples of
(\ref{two-ex}) are reversible in $\GL(3,\ZZ)$. 

In this Section, we are concerned with determining when reversibility
can occur in $\GL(n,\ZZ)$, and what we can say about the nature of the 
reversing symmetry $G \in \GL(n,\ZZ)$, e.g.\ whether it can be taken to be an
involution so that, by Fact \ref{rev-group}, the reversing symmetry group 
${\cal R}(M) \subset \GL(n,\ZZ)$ is a semi-direct product. Note that if
$M \in \GL(n,\RR)$ is reversible, it has been shown that there always exists 
an involutory reversing symmetry \cite[Thm.\ 2.1]{Sev}. Already for 
$\GL(2,\ZZ)$, this is no longer true \cite{Wilson}: the matrix 
$M=\mbox{\tiny $\left(\begin{array}{cc}\! 5 \! & \! 7 \! \\ 
         \! 7 \! & \! 10 \! \end{array}\right)$}$  
is reversible in $\GL(2,\ZZ)$ with the reversing symmetry 
$G=\mbox{\tiny $\left(\begin{array}{cc} \! 0 \! & \! 1 \! \\
         \! -1 \! & \! 0 \! \end{array}\right)$}$
of order 4,
but with no involutory reversing symmetry in $\GL(2,\ZZ)$.

While reversibility was still a frequent phenomenon in $\GL(2,\ZZ)$, 
see \cite{BR97}, it becomes increasingly restrictive with growing $n$.
To see this, let us consider necessary conditions for
reversibility. If $M\in \GL(n,\RR)$ is reversible, then
$M^{-1}=GMG^{-1}$ for some matrix $G$, and $M$ and $M^{-1}$ must
have the same characteristic polynomial, $P(x)$. On the other
hand, the spectrum of $M$ must then be self-reciprocal, i.e.\ with
$\lambda$ also $1/\lambda$ must be an eigenvalue, with matching
multiplicities. Recall that $P(x)=\prod_{i=1}^{n} (x-\lambda_i)$
and observe that
\begin{equation}
\prod_{i=1}^{n} \Bigl(x-\frac{1}{\lambda_i}\Bigr) \; = \;
\frac{(-1)^n x^n}{\det(M)}\,\prod_{i=1}^{n}
\Bigl(\frac{1}{x}-\lambda_i\Bigr).
\end{equation}
But by assumption, $\prod_{i=1}^{n} (x-\lambda_i) =
\prod_{i=1}^{n} (x-\frac{1}{\lambda_i})$, so we arrive at
\begin{prop}\label{rec-poly}
   A necessary condition for the matrix $M\in\GL(n,\RR)$ to be 
   reversible is\/ $\spec(M) = \spec(M^{-1})$ and thus the equation
\begin{equation} \label{reci}
   P(x) \; = \; \frac{(-1)^n x^n}{\det(M)}\,P(1/x)
\end{equation} 
   which we call the self-reciprocity of $P(x)$. \qed
\end{prop}
One immediate consequence is the following.
\begin{coro}\label{det-is-one}
  If $M\in\GL(n,\RR)$ is reversible, $\det(M)=\pm 1$. If, in addition, the
  multiplicity of the eigenvalue $\lambda=-1$ is even $\,($allowing
  for multiplicity $0$ if $-1$ is not an eigenvalue of $M$$\,)$, then
  $\,\det(M)=1$.
\end{coro}
{\sc Proof}: Observe that the reversibility of $M$ implies
$G=MGM$. Taking determinants gives the first assertion because
neither $M$ nor $G$ is singular.

Next, note that $\lambda=\pm 1$ are the only complex numbers with
$\lambda=1/\lambda$. All other eigenvalues come in reciprocal
pairs. Since the determinant is the product over all eigenvalues,
the second statement follows.  \qed

It might be instructive to reformulate Proposition \ref{rec-poly}
and Corollary \ref{det-is-one} in terms of elementary symmetric
polynomials. Let $S_k$, $k=0,1,...\, ,n$, denote the $k$-th elementary 
symmetric polynomial in $n$ indeterminates. The $S_k$ are given by 
$S_0\equiv 1$ and
\begin{equation}
   S_k(x^{}_1,...\, ,x^{}_n) \; = \; 
     \sum_{i^{}_1 < \dots \, < i^{}_k} 
     x^{}_{i^{}_1}\cdot\ldots\cdot x^{}_{i^{}_k}\, .
\end{equation}
They are algebraically independent over $\ZZ$ and have the generating function
\begin{equation}
   \sum_{k=0}^n S_k(x^{}_1,...\, ,x^{}_n) \, t^k \; = \;
   \prod_{i=1}^n \, (1+x^{}_i t)
\end{equation}
where $t$ is another indeterminate.

Now observe that, for a characteristic polynomial $P(x)$, we have 
$$P(x) \; = \; x^n +
  \sum_{k=1}^{n} (-1)^k S_k(\lambda_1,...\, ,\lambda_n)\, x^{n-k}\, .$$
Consequently, $$x^n P(1/x)\; = \; 1+\sum_{k=1}^{n} (-1)^k
S_k(\lambda_1,...\, ,\lambda_n)\, x^k$$ and a comparison with
Proposition \ref{rec-poly} reveals that
\begin{equation}
S_k(\lambda_1,...\, ,\lambda_n)\; = \;\det(M)\cdot S_{n-k}
(\lambda_1,...\, ,\lambda_n)
\end{equation}
for all $0\le k\le n$. For $k=0$, this is just the statement that
$\det(M)=\pm 1$. Note that the elementary symmetric polynomials, when
evaluated at the roots of $P(x)$, reproduce (up to a sign) the entries of
the last row of the left companion matrix (\ref{left-compa}).

Turning now to the reversibility of matrices in $\GL(n,\ZZ)$, we first
observe that, generically, the reversible cases can only occur when $n$ 
is even and $\det(M)=+1$:
\begin{prop}\label{criteria}
  Consider $M\in \GL(n,\ZZ)$ and let $P(x)$ be the characteristic 
  polynomial of $M$. If $n > 1$ is odd or $\det(M)=-1$ we have:
\begin{itemize}
\item[$(a)$] $\;$if $M$ is reversible in $\GL(n,\ZZ)$, $P(x)$ is reducible 
             over $\ZZ$, 
     and the spectrum of $M$ contains $1$ or $-1$;
\item[$(b)$] $\;$if $P(x)$ is irreducible over $\ZZ$, $M$ cannot be 
             reversible in $\GL(n,\ZZ)$.
\end{itemize}
\end{prop}
{\sc Proof}: If $M$ is reversible, $\lambda\in\spec(M)$ implies
$1/\lambda\in\spec(M)$, so the eigenvalues are either $\pm 1$ or
have to come in pairs, $\lambda\neq 1/\lambda$. If $n$ is odd,
we must have at least one eigenvalue that is $\pm 1$, and that
gives a factor $(x\mp 1)$ in $P(x)$.
On the other hand, if $\det(M)=-1$, we must
have at least one eigenvalue that is $-1$ which gives a factor
$(x+1)$ in $P(x)$. In both cases, one notes that whenever $\pm 1$
is a zero of a polynomial over $\ZZ$, factoring out $(x\mp 1)$ can be done
over $\ZZ$.

Conversely, if $P(x)$ is irreducible over $\ZZ$, $\spec(M)$ cannot
contain an eigenvalue of the form $\pm 1$, and $n$ odd or
$\det(M)=-1$ is then
incompatible with $M$ being reversible. \qed

\smallskip
It is clear from this that reversibility is rather restrictive. If
$P(x)$ splits into irreducible components $P_i(x)$, then each is
subject to the constraints described above, or has to be matched
with its reciprocal partner polynomial -- if that would be an
integer polynomial at all. In particular, if $P(x)$ is reducible
but contains an isolated irreducible factor of odd order $\ge 3$, or
of even order with constant term $-1$, reversibility of $M$ is ruled out. 
For example, this confirms that $M_1$ and $M_2$ of (\ref{two-ex}) are not 
reversible (in fact, they are not even reversible in $\GL(3,\RR)$).

The key problem in deciding upon similarity of $M$ and $M^{-1}$ in
$\GL(n,\ZZ)$ is that $\ZZ$ is not a field. But it is clear that the
corresponding similarity within $\GL(n,\QQ)$ (the matrix entries now belonging 
to the field of rationals) is both a necessary condition and a much easier
problem. It would not help to further extend $\QQ$ to $\RR$ due to the
following result, see \cite[Cor.\ XIV.2.3]{Lang}.
\begin{fact} \label{no-help}
  A matrix $M\in\GL(n,\ZZ)$ is similar to $M^{-1}$ within the group 
  $\GL(n,\RR)$ if and only if this is already the case in $\GL(n,\QQ)$. \qed
\end{fact}

In the light of this, let us first recall some facts about normal forms 
over $\QQ$, where similarity is (in theory) a decideable problem. 
The normal form of a matrix $M$ is based on its {\em polynomial invariants\/},
or {\em invariants\/} for short, see \cite[Sec.\ XIV.2]{Lang}.
They are often also called the {\em invariant factors\/} of $M$
(or, more explicitly, of $(x\Id\!-\!M)$), compare 
\cite[Def.\ 4.4.6]{adwein}, meaning certain polynomials that derive
from the matrix $(x\Id-M)$, see below. The following result is a direct 
consequence of \cite[Thm.\ XIV.2.6]{Lang} or \cite[Thm.\ 5.3.3]{adwein}.
\begin{fact} \label{poly-sim}
   Two matrices in ${\rm Mat}(n,\QQ)$ are similar in $\GL(n,\QQ)$ 
   if and only if they have the same polynomial invariants.
   In particular, this applies to $M$ and $M^{-1}$ for
   any $M\in\GL(n,\ZZ)$. \qed
\end{fact}

Let us briefly recall how the polynomial invariants $q^{}_1,...\, ,q^{}_r$ of
a matrix $M\in{\rm Mat}(n,\ZZ)$ can be found, where $r\le n$ is a
uniquely determined integer that depends on $M$. We formulate this for
integer matrices, but it applies, with little change, also to rational ones.
Set $p^{}_0=1$ and let $p^{}_k$ (for $1\le k\le n$) be the greatest common 
divisor of all minors of $(x\Id\!-\!M)$ of order $k$, so that $p^{}_k$ clearly 
divides $p^{}_{k+1}$, and $p^{}_n = P(x) = \det(x\Id\!-\! M)$. Let $\ell$ 
denote the largest integer $k$ for which $p^{}_k=1$ and define 
$q^{}_i = p^{}_{\ell+i}/p^{}_{\ell+i-1}$, where $1\le i\le r=n-\ell$.
These polynomials over $\ZZ$ are the polynomial invariants of $M$ and 
satisfy the following divisibility property:
\begin{equation}  \label{divi}
     q^{}_i \; | \; q^{}_{i+1}\, .
\end{equation}
The prime factors of $q^{}_i$ over $\ZZ$, taken with their multiplicity, are
called its {\em elementary divisors\/}. The product of the invariant
factors of $M$ (equivalently, the product of all their elementary
divisors) gives the characteristic polynomial $P(x)$ of $M$.
Furthermore, the minimum polynomial $Q(x)$ of $M$ is given by
$q^{}_r$, or, equivalently, by the characteristic polynomial $P(x)$
divided by $p^{}_{n-1}$.

Note that a systematic way to find the invariant factors of $M$
is to bring $(x\Id-M)$, seen as a matrix over the principal ideal
domain $\ZZ[x]$, into its so-called {\em Smith normal form\/}, 
see \cite[Ch.\ 5.3]{adwein} for details.
This is a diagonal matrix in of the form
${\rm diag}(1,...\, ,1, q^{}_1(x), ... \, , q^{}_r(x))$.
{}For large $n$, calculating 
this form can be a computationally difficult exercise; for small $n$,
the Smith normal form can be found from algebraic program packages.
Nevertheless, significantly, the invariant factors completely determine the
{\em Frobenius normal form\/} of the matrix $M$:
\begin{fact} \label{Frobenius}
   Let $M\in{\rm Mat}(n,\ZZ)$ have polynomial invariants  
   $q^{}_1,...\, , q^{}_r$ of degrees 
   $n^{}_1,...\, , n^{}_r$, with $n^{}_1+ \ldots +n^{}_r=n$. Then $M$
   is similar, in $\GL(n,\QQ)$, to a block diagonal matrix
   $[B_1,...\, , B_r]$ where $B_i$ is the $n_i \times n_i$
   left companion matrix of the polynomial $q^{}_i$.  \qed
\end{fact}

The existence of a block diagonal matrix similar to $M$ is equivalent to
the statement that $M$ leaves invariant a set of (cyclic) subspaces
of $\QQ^n$ with respective dimensions $n^{}_1,...\, , n^{}_r$, see
\cite[Thm.\ XIV.2.1]{Lang} for details. One can actually give more 
refined normal forms by using the elementary divisors of each invariant 
to replace the diagonal blocks $B_i$ with subblock decompositions
based upon the elementary divisors and their multiplicities.
Combining Fact \ref{poly-sim} and Fact \ref{Frobenius}, matrices with 
the same polynomial invariants can both be brought to the same normal 
form and thus are similar.

The normal form of Fact \ref{Frobenius} highlights the left companion matrices
$B_i$. For what follows, we are interested in the reversibility
of such matrices. In this respect, let $M^{(\ell)}$ and $M^{(r)}$ be
the left and right companion matrices corresponding to a polynomial
$P(x)$. Suppose $P(x)$ conforms to the reciprocity condition (\ref{reci}) 
of Proposition \ref{rec-poly}. Then one can check that $M^{(r)}$ is the
{\em inverse\/} of $M^{(\ell)}$. But we already know from (\ref{comp-conj}) 
that $M^{(r)}=RM^{(\ell)}R^{-1}$ where $R=R^{-1}$ is the involution from
(\ref{R-matrix}).  Combining this with the normal form above, we obtain:
\begin{theorem} \label{Q-rev}
   Let $M\in\GL(n,\ZZ)$. Then, $M$ is reversible in $\GL(n,\QQ)$ if and only 
   if each of the polynomial invariants of $M$ satisfies the reciprocity 
   condition $(\ref{reci})$ separately. In this situation,
   the reversing symmetry can be chosen to be an involution.
\end{theorem}

{\sc Proof}:
By Fact \ref{Frobenius},  $M=SDS^{-1}$ where $S \in \GL(n,\QQ)$ and
$D \in \GL(n,\ZZ)$ is a block diagonal matrix of the form
$D=[B_1,...\, , B_r]$, where $r \ge 1$ and $B_i \in \GL(n_i,\ZZ)$ is the
left companion matrix corresponding to the invariant $q^{}_i$ of degree 
$n_i$. It follows that $M^{-1}=SD^{-1}S^{-1}$ where
$D^{-1}=[B_1^{-1},...\, , B_r^{-1}]$. Consequently, $M$
and $M^{-1}$ are similar if and only if $D$ and $D^{-1}$ are similar.

Suppose that each of the polynomial invariants satisfies the condition
(\ref{reci}). Then, from the remark before Theorem \ref{Q-rev}, $B_i^{-1}$
is the right companion matrix corresponding to $q^{}_i$ and is
similar to $B_i$ via the involution $R_i \in \GL(n_i,\ZZ)$ which
consists of $1$'s on its anti-diagonal, as in (\ref{R-matrix}). 
It follows that $D$ is similar to $D^{-1}$ via the block diagonal involution
$R:=[R^{}_1,...\, , R^{}_r]$, and so $M$ and $M^{-1}$ are similar by the
involution $S R S^{-1}$.

On the other hand, suppose that $M$ is similar to $M^{-1}$
in $\GL(n,\QQ)$. Hence, the corresponding block diagonal matrices $D$
and $D^{-1}$ are also similar in $\GL(n,\QQ)$, via some
element $G$. Now, to each block $B_i$ of $D$ corresponds an invariant
vector subspace $V_i$ of $\QQ^n$ of dimension $n_i$. 
A subspace $V_i$ is thus either mapped
by $G$ to itself (it is a {\em symmetric\/} subspace) or
to another subspace $V_j $ of the same dimension. In
the first case, $B_i$ must be conjugate to its inverse
via the restriction of $G$ to $V_i$. This means that
the characteristic polynomial of $B_i$, which is $q^{}_i$,
must satisfy the condition (\ref{reci}). On the other hand,
if $V_i$ is mapped to $V_j $ by $G$ with $n_i=n_j $, it follows
that the invariants $q^{}_i$ and $q^{}_j $ differ by at
most a sign. They thus share the same eigenvalues and must
each satisfy condition (\ref{reci}) on their own.  \qed

\smallskip
If $M$ has only one non-trivial invariant, $q^{}_1(x)$, it follows from 
the above discussion that its characteristic polynomial $P(x)$
coincides with its minimal polynomial $Q(x)$ and both equal
$q^{}_1(x)$ (so $M$ is cyclic from Fact \ref{facts}). Conversely, $M$
cyclic means it has only one invariant. The previous Theorem now gives:
\begin{coro} \label{coro5}
   If $M\in\GL(n,\ZZ)$ has only one polynomial invariant, in
   particular if the characteristic polynomial
   $P(x)$ is irreducible over\/ $\ZZ$, then
   $M$ is reversible in $\GL(n,\QQ)$ if and only if its
   characteristic polynomial $P(x)$ satisfies the
   reciprocity condition $(\ref{reci})$.
   \qed
\end{coro}

Note that Theorem \ref{Q-rev} and Corollary \ref{coro5} do not extend 
to requiring, for reversible $M$, that the elementary divisors within an 
invariant polynomial satisfy (\ref{reci}).
For example, any matrix in $\GL(4,\ZZ)$ with one invariant polynomial 
$q^{}_1(x)=P(x)=(x^2-x-1)(x^2+x-1)$ is reversible in $\GL(4,\QQ)$ 
although the elementary divisors separately violate (\ref{reci}).

Let us give some illustrations of the use of Theorem \ref{Q-rev} and
Corollary \ref{coro5} for small values of $n$. These results 
show that all $M \in {\rm SL}(2,\ZZ)$
are reversible in $\GL(2,\QQ)$ because their invariant
factors fall into one of the following cases:
\begin{enumerate}
\item $q^{}_1(x)=q^{}_2(x)=(x\pm 1)$, so $r=2$ and $M = \mp\Id$; 
\item $q^{}_1(x)=P(x)=x^2 - {\rm tr}(M) x + 1$, so $r=1$ and $P(x)$
      is self-reciprocal.
\end{enumerate}
Yet, we know from \cite{BR97} that ${\rm SL}(2,\ZZ)$ matrices exist that 
are not reversible in $\GL(2,\ZZ)$ (see also Section \ref{gen-symm} below
for further discussion). Furthermore, if $M \in \GL(2,\ZZ)$ with
$\det M = -1$, then it can only have one polynomial invariant,
$q^{}_1(x)=P(x)=x^2- {\rm tr}(M) x -1$. 
By Proposition \ref{criteria} or Corollary \ref{coro5}, $M$ is reversible 
in $\GL(n,\ZZ)$ if and only if $M$ has eigenvalues $\lambda = \pm 1$ and
$q_1(x)$ factors into $(x-1)(x+1)$, in
which case $M\in \GL(2,\ZZ)$ is an involution and reversible,
with reversing symmetry as itself. This approach gives another
way of retrieving some of the results of \cite{BR97} on the reversibility
in $\GL(2,\ZZ)$.

Turning to $\GL(3,\ZZ)$, Proposition \ref{criteria} implies that if $M$
is reversible then $P(x)$ must have a factor $(x \pm 1)$.
In other words, $(x \pm 1)^i$, for some $1 \le i \le 3$, is an elementary 
divisor of $P(x)$.
If $M$ has more than one polynomial invariant, then the divisibility
property (\ref{divi}) implies that $P(x)$ completely decomposes into
a product of 3 factors of the form $(x \pm 1)$. Generically,
however, this will not happen and instead  $M$ is reversible
in $\GL(3,\QQ)$ if and only if it has only one invariant of the
form $q^{}_1(x)=P(x)=(x \pm 1)(x^2-({\rm tr}(M) \pm 1) x +1)$.

{}For $\GL(4,\ZZ)$, elements with three or four polynomial invariants
are diagonal matrices with $+1$'s or (an even number of)
$-1$'s on the diagonal. They are all reversible. Reversible elements
with two invariants must have $q^{}_1(x)=( x \pm 1)$
and $q^{}_2(x)= (x \pm 1) (x^2 - (\tr(M) \pm 2) x +1)$, or
$q^{}_1(x)=q^{}_2(x)$, a monic quadratic with constant term $+1$.

As $n$ increases, Theorem \ref{Q-rev} can exclude many matrices from being
reversible in $\GL(n,\ZZ)$ because they are not reversible in $\GL(n,\QQ)$.  
In particular, we can ask for an example $M$ with more than one
polynomial invariant where the characteristic polynomial
$P(x)$ satisfies (\ref{reci}), yet $M$ is irreversible in $\GL(n,\QQ)$
because it violates Theorem \ref{Q-rev}. If we take $n \ge 2$
and even, the first possibility appears in $\GL(8,\ZZ)$. For
example, we can take a matrix with invariants
$q^{}_1(x)=x^2-x-1$ and $q^{}_2(x)=(x^2-x-1)(x^2+x-1)^2$. Both
polynomials violate the reciprocity condition (\ref{reci}) but
they compensate each other so that their product, the characteristic
polynomial, does satisfy the condition. The Frobenius normal form
with these invariants is the block diagonal matrix $M=[B_1,B_2]$,
\begin{equation} \label{8-cat}
M \; = \; \begin{pmatrix}
               0 & 1 & 0 & 0 & 0 & 0 & 0 & 0 \\
               1 & 1 & 0 & 0 & 0 & 0 & 0 & 0 \\
               0 & 0 & 0 & 1 & 0 & 0 & 0 & 0 \\
               0 & 0 & 0 & 0 & 1 & 0 & 0 & 0 \\
	       0 & 0 & 0 & 0 & 0 & 1 & 0 & 0 \\
               0 & 0 & 0 & 0 & 0 & 0 & 1 & 0 \\
               0 & 0 & 0 & 0 & 0 & 0 & 0 & 1 \\
               0 & $\hphantom{-}$0 & $\hphantom{-}$1 
                 &$-$1 &$-$4 & $\hphantom{-}$3 & $\hphantom{-}$4 &$-$1
               \end{pmatrix}.
\end{equation}
This matrix is not reversible in $\GL(8,\ZZ)$ (although its square is,
see Section \ref{projective} below).

The remaining problem is now to find possible reversibility
in $\GL(n,\ZZ)$ of unimodular matrices that are already reversible
in $\GL(n,\QQ)$. For $n=2$, we were able to solve this problem
using a special algebraic structure (the amalgamated free product)
of ${\rm PGL}(2,\ZZ)$. This structure is not available for $n \ge 3$.

Important examples of unimodular integer matrices which are reversible
in $\GL(n,\QQ)$ are the {\em symplectic\/} matrices
${\rm Sp}(2n,\ZZ) \subset {\rm SL}(2n,\ZZ)$. Recall that a symplectic matrix
$M \in {\rm Sp}(2n,\RR)$ satisfies $M^t J M = J$ where $M^t$ denotes the
transpose of $M$ and $J$ is the $2n\! \times\! 2n$ integer block matrix 
$$ J \; = \; \begin{pmatrix} \bs{0} & \Id \\ $-$\Id & \bs{0} \end{pmatrix} $$
of order $4$. Since in general
\begin{equation} \label{symplectic}
   M^t J M \; = \; J \quad \Rightarrow \quad
   M^t \; = \; J M^{-1} J^{-1}\, ,
\end{equation}
it follows that $M \in {\rm Sp}(2n,\RR)$ is reversible if and only if $M$ is
similar to $M^t$ in $\GL(2n,\RR)$. But any invertible square matrix with
entries in a field $F$ is similar to its transpose in $\GL(m,F)$, see
\cite[Prop.\ 5.3.7]{adwein} (but this need not be true e.g.\ in
$\GL(m,\ZZ)$). In particular,
$M \in {\rm Sp}(2n,\ZZ)$ is reversible in $\GL(2n,\QQ)$ and its invariant
factors will all satisfy (\ref{reci}). Also, it is clear from 
(\ref{symplectic}) that if $M$ is
symplectic and symmetric, then  $M^t=M= J M^{-1} J^{-1}$ so that $M$ is
actually reversible in $\GL(2n,\ZZ)$. 

Orthogonal integer matrices
$U \in \GL(n,\ZZ)$ which satisfy $U U^t=\Id$ are other examples of unimodular
integer matrices which are always reversible in $\GL(n,\QQ)$ 
(since $U^{-1}=U^t$ and, as before, $U^t$ and $U$ are similar in 
$\GL(n,\QQ)$).

We make some further remarks on the problem of deciding reversibility
within $\GL(n,\ZZ)$.
Firstly, note that in the proof of Theorem \ref{Q-rev} above, when we used 
the reversibility of the left companion matrices with characteristic polynomials
satisfying (\ref{reci}), this reversibility was in $\GL(n,\ZZ)$ itself. 
Furthermore, the reversing symmetry $R$ was an involution, so by 
Fact~\ref{rev-group} we have the following result.
\begin{theorem}
  {}For each integer polynomial $P(x)$ of degree $n$ that 
  satisfies the necessary 
  self-reciprocity condition $(\ref{reci})$ for reversibility, there is at 
  least one reversible matrix class in $\GL(n,\ZZ)$, represented by the 
  left companion matrix $M^{(\ell)}$ of $P(x)$, and we have 
  ${\mathcal R}(M^{(\ell)}) = {\mathcal S}(M^{(\ell)}) \times_s C_2$.  \qed
\end{theorem}

Secondly, by Fact \ref{classnumber},
the number of representing matrix classes of an irreducible
characteristic polynomial $P(x)$ equals the class number of the
order $\ZZ[\alpha]$, with $\alpha$ any of the roots of $P(x)$. If
this class number is one, there is {\em only\/} the class
represented by the companion matrix. If the class number is two,
one is the companion matrix class which we know to be reversible
if the spectrum is self-reciprocal. But then, the other class must
also be reversible because there is no further partner left.

Class numbers are widely studied in algebraic number theory, and
one can find both extensive tables in books (e.g.\ see \cite{Hasse,PZ}) 
and also various program packages to calculate them, e.g.\ the program
package KANT\footnote{See {\sf http://www.math.TU-Berlin.de/$\sim$kant/}.}.

Let us add another example, of rather different flavour, and look at
an interesting class of algebraic
integers, the so-called {\em Salem numbers\/}. They are the
algebraic integers $\alpha>1$ with all conjugates $\alpha'$ having
modulus $|\alpha'|\le 1$ and with at least one conjugate on the
unit circle, see \cite[Ch.\ III.3]{Salem} for details. 
Salem's Theorem then says that their
degree is always even, that $\alpha$ and $1/\alpha$ are the only real
conjugates, and that all other conjugates are on the unit circle.
In particular, a Salem number is a unit.
Putting our above results to work, we get
\begin{coro}\label{salem}
  Each Salem number occurs as the eigenvalue of a reversible toral
  automorphism. If $\alpha$ is a Salem number of degree $n=2m$, and
  $M$ a corresponding $\GL(n,\ZZ)$ matrix, then 
  $\,{\mathcal S}(M)\simeq C_2\times\ZZ^m$. \hfill $\square$
\end{coro}

The polynomial $P(x)=x^4-2x^3-2x^2-2x+1$ provides one of the
simplest examples. Its roots are $\tau\pm\sqrt{\tau}$ (both real)
and $-(\tau-1)\pm i\sqrt{\tau-1}$ (both on the unit circle), where
$\tau=\bigl(1+\sqrt{5}\,\bigr)/2$ is the golden number. The
corresponding companion matrix $M$ is not of finite order, and the
reversing symmetry group is thus
${\mathcal R}(M)\simeq (C_2\times \ZZ^2)\times_s C_2$.

Let us come back to the general discussion and ask for the properties
of reversing symmetries. Let $G$ be a reversing symmetry of $M$, so
that $GM = M^{-1}G$ and hence also $GM^n = M^{-n}G$ for all $n\in\ZZ$.
If $p(x)$ is any polynomial, we then also get $G p(M) = p(M^{-1}) G$.
Since $G$ is a reversing symmetry, $G^2$ is a symmetry. If we now assume 
that $M\in\GL(n,\ZZ)$ has minimal polynomial of degree $n$, we get
$G^2 \in \QQ[M]$ from Fact \ref{facts}, i.e.\ $G^2 = q(M)$ for some
$q$ with coefficients in $\QQ$. Consequently,
$G q(M) G^{-1} = q(M^{-1})$, but also $G q(M) G^{-1} = G^2 = q(M)$,
so that $q(M) = q(M^{-1})$.

If $q(M)$ is a monomial, i.e.\ $q(M)=M^{\ell}$ for some $\ell$, then
$M^{2\ell}=1$ and hence $G^4=1$. This case is also discussed in
\cite[Prop.\ 2(i)]{G}. It clearly extends to the situation that 
$q(M^{-1}) = \bigl(q(M)\bigr)^{-1}$, which is more general.
Apart from this, we recall the following result from \cite[p.\ 21]{G}
(its proof, which was only contained in the preprint version of
\cite{G}, is a coset counting argument).
\begin{fact} \label{finite-order}
   Let $M$ be of infinite order. If the factor group 
   ${\mathcal S}(M)/\langle M\rangle$ is finite,
   then any reversing symmetry $G$ of $M$ must be of finite order,
   and $G^{2k}=\Id$ for some integer $k$ that divides the order
   of the factor group. \qed
\end{fact}
Let us only add that, in line with our above argument, one first
obtains $G^{2k}\in\langle M\rangle$ and hence $G^{4k}=\Id$. But
$\langle M\rangle\simeq \ZZ$ by assumption, so it cannot have a
subgroup of order 2, and thus already $G^{2k}=\Id$.

Fact \ref{finite-order} certainly applies to our scenario whenever $M$ 
is not of finite order, but ${\mathcal S}(M)$ Abelian and or rank 1.
This type of result is helpful because it restricts the search for reversing 
symmetries to one among elements of finite order. It is certainly possible to
extend the result to other cases, but in general one has to expect
reversing symmetries of infinite order, in particular if the rank of
${\mathcal S}(M)$ is $\ge 2$. Even then some results are possible
because it would be sufficient to know whether reversibility implied
the existence of {\em some\/} reversing symmetries of finite order.
However, this question is more involved and thus postponed.

\section{Extensions and further directions}\label{further}

In this Section, we will summarize some additional aspects of our
analysis, namely the extension of symmetries to affine mappings,
the modifications needed to treat the related situation of the
projective matrix group $\PGL(n,\ZZ)$, and the extension of
(reversing) symmetries from a group setting to that of (matrix) 
rings or semi-groups.

\subsection{Extension to affine transformations}

So far, we have mainly discussed linear transformations (w.r.t.\ the torus), 
but it is an interesting question what happens if one extends the search 
for (reversing)
symmetries to the group of {\em affine\/} transformations. Since both
arguments and results are the exact analogues of those for the
case $n=2$ as derived in \cite{BR97}, we will be very brief here.

In Euclidean $n$-space, the group of affine transformations is the
semi-direct product ${\mathcal G}_{a} = \RR^n \times_s
\GL(n,\RR)$, with $\RR^n$ being the normal subgroup. Elements are
written as $(t,M)$ with $t \in \RR^n$ and $M \in \GL(n,\RR)$, and
the product of two transformations is $(t,M) \cdot (t',M') = (t +
M t', M M')$. The neutral element is $(0,\Id\,)$, and we have
$(t,M)^{-1} = (- M^{-1}t,M^{-1})$.

If we now observe that $\TT^n = \RR^n/\ZZ^n$, it is immediately clear that
the affine transformations of $\TT^n$ form the group
\begin{equation} \label{aff5}
    {\mathcal G}_{a}^{\TT^n} \; = \; \TT^n \times_s \GL(n,\ZZ)
\end{equation}
which is still a semi-direct product. Here, $\TT^n$ can be written as
$[0,1)^n$ with addition mod.\ 1, and the product of transformations is
modified accordingly.

If we now ask for an affine (reversing) symmetry of a matrix $M$
(now being identified with the element $(0,M) \in {\mathcal
G}_{a}^{\TT^n}$) we find
\begin{prop} \label{affine-prop}
   The affine transformation $(t,G)$ is a $($reversing$\/)$
   symmetry of the toral automorphism $(0,M)$
   if and only if
\begin{itemize}
   \item[$(a)$] $\;G$ is a $($reversing$\,)$ symmetry of $M$ in $\GL(n,\ZZ)$ and
   \item[$(b)$] $\;Mt=t$ $( \mbox{mod } 1)$.
\end{itemize}
\end{prop}
{\sc Proof}: We have $(t,G) \cdot (0,M) = (t,GM)$ and also
$(0,M^{\pm 1}) \cdot (t,G) = (M^{\pm 1}t,M^{\pm 1}G)$. But then,
the statement follows from the uniqueness of factorization in
semi-direct products. \qed

{}From the condition $Mt=t$ (mod 1) it is clear that we need not consider all
translations in $\TT^n$ but only those with {\em rational\/} components, 
which we denote as $\Lambda_{\infty}$. For many concrete problems, it would 
actually be even more appropriate to restrict
to discrete sublattices, e.g.\ to the so-called $q$-division points
$\Lambda_q \simeq (C_q)^n$ which consists of all rational points with 
denominator $q$. We will not follow this idea here, however.

{}From the above result, it is clear that we can get (reversing)
$k$-symmetries (recall the definitions from Section \ref{sec-ksym}). 
In fact, the equation $M^k t = t$ on the torus has
$a_k = |\det(M^k - \Id)|$ different solutions provided no
eigenvalue of $M^k$ is 1. Clearly,
\begin{equation} \label{aff9}
    a_k^{} \; = \; \sum_{\ell | k} \ell \cdot c_{\ell}^{}
\end{equation}
where $c_{\ell}^{}$ counts the true orbits of length $\ell$,
and the M\"obius inversion formula gives
\begin{equation} \label{aff10}
    c_k^{} \; = \; \frac{1}{k} \sum_{\ell | k}
              \mu (\mbox{\small $\frac{k}{\ell}$}) \cdot a_{\ell}
\end{equation}
with the M\"obius function $\mu(m)$ \cite[p.\ 29]{Cohn}.
If $c_k^{}$ is positive for some $k$, we get a $k$-symmetry
(and, hence, eventually a reversing $k$-symmetry) of $M$.
These numbers can easily be calculated explicitly, where a very natural
tool is provided by the so-called dynamical or Artin-Mazur
$\zeta$-functions \cite{EI}. Here, the $a_k^{}$'s can be
extracted from the series expansion of the logarithm of the
$\zeta$-function, while the $c_k^{}$'s appear as exponents
of the factors of the Euler product expansion of the
$\zeta$-function itself.

\subsection{The case of $\PGL(n,\ZZ)$}  \label{projective}

Let us start by the observation that $\PGL(n,\ZZ)$ can be
described via quotienting w.r.t.~$\{\pm \, \Id\}$, i.e.\
$$\PGL(n,\ZZ) \simeq \GL(n,\ZZ)/\{\pm\, \Id\}\, .$$
 In other words, rather than consider single matrices $M$,
one has to consider pairs, $[M]:=\{\pm M\}$. Let us write
${\mathcal S}[M]$ for the new $\PGL$ case and keep the old
notation for the $\GL$ situation treated above.

The modification needed for the symmetry analysis given above is
then actually fairly trivial, as we always had $\pm \Id$ among
them, and we can simply factor that out. So, we get
\begin{equation}
   {\mathcal S}[M] \; \simeq \; {\mathcal S}(M)/\{\pm \Id\}\,.
\end{equation}

The case of reversing symmetries, however, requires some care.
Since we now calculate mod $\pm \Id$, a projective matrix $[M]$
can also be reversible through $GMG^{-1}=-M^{-1}$. But if this
happens, the square, $M^2$, is again reversible in the old sense.
This mechanism can (and will) give rise to reversing 2-symmetries
(recall Section \ref{sec-ksym} for the definition) in $\GL(n,\ZZ)$.
Whereas \cite{BR97} gave examples in $\GL(2,\ZZ)$, 
Eq.~(\ref{symplectic}) shows that
skew-symmetric symplectic matrices satisfying $M^t=-M$
are not reversible in $\GL(2n,\ZZ)$ whereas their squares are reversible.
Also, the example (\ref{8-cat}) is reversible in ${\rm PGL}(8,\ZZ)$ since 
a $G \in \GL(8,\ZZ)$ can be found, by direct calculation, that satisfies
the relation $GMG^{-1}=-M^{-1}$.

Let us finally check what happens in the extension to affine
transformations. In complete analogy to the case $n=2$, see
\cite{BR97}, one can show that the corresponding affine group is
the semidirect product $\Lambda_2 \times_s {\rm PGL}(n,\ZZ)$ with
$\Lambda_2$ the 2-division points. This really is the consequence
of identifying $x$ with $-x$ on $\TT^n$, and $\Lambda_2$ is the
set of $2^n$ translations that satisfy the condition $t  =  -t
\quad \mbox{(mod 1)} $.

Now, the above Proposition \ref{affine-prop} applies to the case of 
$\PGL$-matrices in very much the same way, just the possible translations 
$t$ are restricted to the 2-division points.

\subsection{Symmetries among general integer matrices} \label{gen-symm}

{}For most of this article, we have focused on matrices in
$\GL(n,\ZZ)$ and their symmetries within the same group. However,
none of the proofs given above depends on that restriction, and one 
can indeed also treat the case that both $M$ and its symmetries are
allowed to live in the larger set ${\rm Mat}(n,\ZZ)$ which is no
longer a group w.r.t.\ multiplication, but a ring.  Nevertheless,
we will continue to use the symbol ${\mathcal S}(M)$, now meaning
$${\mathcal S}(M) \; := \;
    \{ G \in {\rm Mat}(n,\ZZ)\mid [M,G]=0 \}\, .$$
The most obvious extended symmetries which one gets in ${\rm Mat}(n,\ZZ)$
are the integer multiples of the identity, but there really
is a hierarchy of objects to look at, and it is most transparent
if one phrases the situation for a matrix in ${\rm Mat}(n,\QQ)$ first:
\begin{fact} \label{rational-facts}
  Let $M\in{\rm Mat}(n,\QQ)$ and let its characteristic polynomial $P(x)$
  be irreducible over $\QQ$. Let $\lambda$ be any of its roots and
  $K=\QQ(\lambda)$ the corresponding algebraic number field.
  Then the following statements hold.
\begin{itemize}
\item[$(a)$] $\;{\rm cent}^{}_{{\rm Mat}(n,\QQ)} (M) \; \simeq \; K$.
\item[$(b)$] $\;{\rm cent}^{}_{\GL(n,\QQ)} (M) \;\; \simeq \; K^*$, where 
             $K^* = K\setminus\{0\}$.
\item[$(c)$] $\;{\rm cent}^{}_{{\rm Mat}(n,\ZZ)} (M) \; \simeq \; \OO$, 
             where $\OO$ is an order in $K$.
\item[$(d)$] $\;{\rm cent}^{}_{\GL(n,\ZZ)} (M) \;\; \simeq \; \OO^{\times}$. 
            \qed
\end{itemize}
\end{fact}
The proof is a slight variation of what we did for Fact \ref{facts}
and Lemma \ref{mainlemma}, and need not be spelled out again. It is
important to note that the order $\OO$ appearing here, as mentioned
before, in general is {\em not\/} the maximal order of $K$, though it
contains $\ZZ[\lambda]$. The following is now an immediate consequence.

\begin{prop} \label{gen-order}
  Let $M$ be an integer matrix with irreducible characteristic
  polynomial $P(x)$. Let $\lambda$ be a root of $P(x)$, and let
  $\OOmax$ be the maximal order in $\QQ(\lambda)$. Then,
  ${\mathcal S}(M)$ is both a $\ZZ$-module and a ring,  and
  isomorphic to an order $\OO$ that satisfies\/ $\ZZ[\lambda]
  \subset \OO \subset \OOmax$. \qed
\end{prop}
Note that there is now also a natural extension to the case of simple 
matrices $M$, compare Theorem \ref{symm-thm1} and Eq.~(\ref{gen-ring}),
but we omit further details here. Also, from Fact \ref{facts} it is clear
that ${\mathcal S}(M)$ is Abelian if and only if the minimal polynomial of
$M$ has degree $n$.

As to reversibility, this new point of view requires some thought.
By Corollary \ref{det-is-one}, $M$ reversible implies $\det(M)=\pm
1$, so reversible integer matrices are restricted to $\GL(n,\ZZ)$.
It would then not be unnatural to also insist on the existence of
at least one unimodular matrix $G$ with $M^{-1}=GMG^{-1}$, and
reversibility is basically as above, except that, if we enlarge
the symmetries of $M$ from subgroups of $\GL(n,\ZZ)$ to subrings of
${\rm Mat}(n,\ZZ)$, reversing symmetries get enlarged accordingly.
Note, however, that the ring structure is lost: the sum of a
symmetry and a reversing symmetry is not a meaningful operation in
this context. Together, they only form a monoid, i.e.\ a
semi-group with unit element.

To go one step further, one could then also rewrite the reversibility
condition as $$G \; = \; M G M$$ and only demand that $G$ is
non-singular, to avoid pathologies with projections to subspaces
and to keep the statement of Corollary \ref{det-is-one}. Note that
this is a slightly weaker form of reversibility, as one does not
assume that $G^{-1}$ is a meaningful mapping in this context. In
particular, it is clear that there is then no need any more to
restrict $G$ to unimodular integer matrices, so that now $G\in\GL(n,\QQ)$.
This will, in general, lead to new cases of (weak) reversibility, as is 
apparent from the explicit constructions in \cite{BR97}.

To give a concrete example, consider the matrices
\begin{equation}
   M \; = \; \begin{pmatrix} 4 & 9 \\ 7 & 16 \end{pmatrix}
   \qquad \mbox{and} \qquad
   G \; = \; \begin{pmatrix} 3 & \, 0 \\ 4 & -3  \end{pmatrix}.
\end{equation}
Then $M$ is the matrix from \cite[Ex.\ 2]{BR97} that was shown to
be irreversible in $\GL(2,\ZZ)$. In fact, the automorphism of
$\TT^2$ induced by $M$ is irreversible even in the larger group of
homeomorphisms of the 2-torus, see \cite{Adler} and \cite[p.\ 9]{AB}
for details on the connection between general and linear homeomorphisms.
Nevertheless, one can
check that $G=MGM$, where $\det(G)=-9$. In other words, $M$ is
reversible in $\GL(2,\QQ)$, as are {\em all\/} elements of ${\rm SL}(2,\ZZ)$
by our previous discussion following Corollary \ref{coro5}.

Note that $G$ does {\em not\/} 
induce a homeomorphism of the 2-torus because $G^{-1}$ is not
an integer matrix. However, $G$ does induce an automorphism on any
lattice of the torus of the form
$$ \Lambda_q \; := \; \{ \, (\mbox{$\frac{m}{q}$}, 
   \mbox{$\frac{n}{q}$})^t\mid
   0 \le m,n < q \}$$
for which $\det(G)\neq 0$ (mod $q$).
This is relevant as a recent study \cite{KM} shows: 
the quantum map which corresponds to $M$ (which, in turn,
corresponds to the action of $M$ on a (Wigner) lattice of the torus)
showed an eigenvalue statistics according to the circular orthogonal 
ensemble (COE) rather than the unitary one (CUE). 
So, even though $M$ does not have a reversing symmetry in the sense of 
Section \ref{sec-rev}, the presence of
``pseudo-symmetries'' such as $G$ still leave their mark!

\subsection*{Acknowledgements}

M.~B.\ would like to thank Peter A.~B.~Pleasants and Alfred Weiss
for several clarifying discussions, Gabriele Nebe for helpful
advice on the literature, and Wilhelm Plesken for suggesting a number 
of improvements. This work was supported by the German
and Australian Research Councils (DFG and ARC).

\end{document}